%menages.tex:Automatic Enumeration of Generalized M\'enage Numbers
%%a Plain TeX file by Doron Zeilberger (x pages)
%begin macros

\baselineskip=14pt
\parskip=10pt

\font\eightrm=cmr8 
\font\eighttt=cmtt8
\magnification=\magstephalf

\def\1{{\overline{1}}}
\def\2{{\overline{2}}}
\parindent=0pt
\overfullrule=0in

\def\frac#1#2{{#1 \over #2}}
%\headline={\rm  \ifodd\pageno  \RightHead  \else  \LeftHead  \fi}
%\def\RightHead{\centerline{
%Title
%}}
%\def\LeftHead{ \centerline{Doron Zeilberger}}
%end macros
\bf
\centerline
{
Automatic Enumeration of Generalized M\'enage Numbers
}
\rm
\bigskip
\centerline{ {\it
Doron 
ZEILBERGER}\footnote{$^1$}
{\eightrm  \raggedright
Department of Mathematics, Rutgers University (New Brunswick),
Hill Center-Busch Campus, 110 Frelinghuysen Rd., Piscataway,
NJ 08854-8019, USA.
%\break
{\eighttt zeilberg  at math dot rutgers dot edu} ,
\hfill \break
{\eighttt http://www.math.rutgers.edu/\~{}zeilberg/} .
Jan. 5, 2014.
Accompanied by Maple package {\eighttt MENAGES}
downloadable from {\eighttt http://www.math.rutgers.edu/\~{}zeilberg/tokhniot/MENAGES} .
Supported in part by the NSF.
}
}
\bigskip

{\it 
\qquad \qquad\qquad 
In fond memory of Alain Lascoux (1944-2013), one of the most CREATIVE and ORIGINAL and INTERESTING
mathematicians that I have ever known}

\bigskip

{\bf One of Alain Lascoux's Favorite Identities}: $0= 1+ (-1)$

Alain Lascoux once said that ``adding and subtracting one'' is one of the most useful tricks in
mathematics, and he was right!

Suppose that you have a large set of {\it guys}, that is easy to count, and there is
a smaller set of {\it good guys} that is hard to count. In other words, the sum
$$
A:=\sum_{ g \in AllGuys} 1 \quad  ,
$$
is easy, but
$$
G:=\sum_{ g \in GoodGuys} 1 \quad  ,
$$
is hard.

Another deep identity in mathematics is, ( for $i$ a non-negative integer)
$$
0^i=0 \,\,\, if \,\,\, i>0 \quad \quad , \quad  \quad 0^0=1 \quad .
$$

If some guy, let's call him Mr. $g$,  is {\bf not} a good guy, it means that he is a {\bf bad guy}, which means that
his set of {\bf sins}, $Sins(g)$,  is {\bf non-empty}, so we have
$$
G=\sum_{g \in AllGuys} 0^{|Sins(g)|} \quad  ,
$$
(the number of elements of a set $S$ is denoted by $|S|$). Now we use Lascoux's favorite identity $0=1+(-1)$ and get
$$
G=\sum_{g \in AllGuys} (1+(-1))^{|Sins(g)|} \quad  .
$$
We need yet another deep identity. For any set $S$, we have:
\footnote{$^2$}
{\eightrm  \raggedright
$$
(1+(-1))^{|S|}= \prod_{s \in S} (1+(-1)) = \sum_{T \subset S} \prod_{s \in S} 1^{\chi (s \not \in T)}  (-1)^{\chi (s \in T)} 
$$
$$
=\sum_{T \subset S} \prod_{s \in S} (-1)^{\chi (s \in T)} =\sum_{T \subset S} \prod_{t \in T} (-1) 
=\sum_{T \subset S} (-1)^{|T|} \quad .
$$
}
$$
(1+(-1))^{|S|}= \sum_{T \subset S} (-1)^{|T|} \quad .
$$
Going back to the formula for $G$ we have
$$
G\, = \,
\sum_{g \in AllGuys} 0^{|Sins(g)|} =\sum_{g \in AllGuys} (1-1)^{|Sins(g)|} 
=\sum_{g \in AllGuys} \sum_{T \subset Sins(g)} (-1)^{|T|} \quad .
$$
By {\bf Discrete Fubini} (aka {\it changing the order of summation}), we have

$$
G \, = \, \sum_{T \subset AllSins}  (-1)^{|T|} \left ( \sum_{\{ g; Sins(g) \supset T \} }  1 \right ) \quad .
$$

Abbreviating, for any set of sins $T$,
$$
A_T \, := \, \sum_{ \{ g; Sins(g) \supset T \} } 1 \, = \, |\{ g \in AllGuys \, ; \, T \subset Sins(g)  \}| \quad ,
$$
the number of guys that are guilty of all the sins in $T$ (but very possibly of other ones as well!),
we have a form of the famous {\bf Inclusion-Exclusion Principle}:
$$
G \, =  \,
\, \sum_{T \subset AllSins}  (-1)^{|T|} A_T \quad .
\eqno(PIE)
$$

{\bf From the Mouth of our Patron Saint: Gian-Carlo Rota}

The {\it patron saint} of enumerative combinatorics,  Gian-Carlo Rota (1932-1999),
started   his {\it seminal} article [Ro] as follows.

``One of the most useful principles of enumeration in discrete probability and combinatorial theory
is the celebrated {\it principle of inclusion-exclusion}. When skillfully applied, this principle
has yielded the solution of many a combinatorial problem $\dots$ '' .

He then goes on to say:

``One frequently notices, however, a wide gap between the bare statement of the principle and
the skill required in recognizing that it applies to a particular problem $\dots$ ''

And indeed, it required lots of human {\it ingenuity} to apply it to many separate problems,
and one of the  most impressive applications was to the enumeration of {\it permutations with restricted positions},
with its accompanying beautiful theory of {\it rook polynomials}, described at length
in the classic [Ri1] of yet-another-guru, John Riordan, and from a more `advanced' viewpoint
in the later classic [St], by still-another guru, Richard Stanley, who used the
{\it transfer-matrix method} (but he should have  mentioned the earlier classic [Ri1]!).

{\bf Back to Inclusion-Exclusion}

In many applications, including the ones discussed in this article, to the counting of {\it discordant permutations},
and more generally, {\it permutations with restricted positions},
it so happens that, in Eq. $(PIE)$,  for many $T$'s, $A_T=0$, since the sins in $T$ are {\it incompatible}, but for those
$T$ for which $A_T \neq 0$, $A_T$ {\bf only} depends on the number of elements of $T$, $|T|$, so we have, in 
these cases
$$
A_T=0 \quad or \quad A_T=f(|T|) \quad ,
$$
for some discrete function $f(k)$.

So things boil down to the {\bf weighted counting} of {\it compatible sets of sins}, $T$, i.e.
sets of sins for which $A_T \neq 0$, according to the weight $(-1)^{k} f(k)$, where $k:=|T|$.

{\bf Enumerating Permutations with Restricted Positions NUMERICALLY} 

Fix a positive integer $n$, say $n=5$. For each $1 \leq i \leq n$, let
there be subsets $A_i$ of $\{1 , \dots , n\}$. 
For example,
$$
A_1 = \{1,4,5\} \quad , \quad 
A_2 = \{3,4 \} \quad , \quad 
A_3 = \{1,3,4 \} \quad , \quad 
A_4 = \{2,3,5 \} \quad , \quad 
A_5 = \{1,4,5 \} \quad .
$$
We want to count the number of $n$-permutations $\pi$ such that
$$
\pi(i) \not \in A_i \quad, \quad i=1 \dots n \quad .
$$

Before  reviewing how to do it via the Principle of Inclusion-Exclusion (Eq. $(PIE)$), and Riordan's lovely {\it rook polynomials},
let's recall the   straightforward, ``{\it positive}'' approach, via {\it permanents}.

{\bf The Positive (Naive) Approach}

Let's define the {\bf complement} sets
$$
B_i := \{ 1, \dots , n \} \backslash A_i \quad ,
$$
so in the above example
$$
B_1 = \{2,3\} \quad , \quad 
B_2 = \{1,2, 5 \} \quad , \quad 
B_3 = \{2,5\} \quad , \quad 
B_4 = \{1,4 \} \quad , \quad 
B_5 = \{2,3 \} \quad .
$$
So our problem is to enumerate the $n$-permutations such that
$$
\pi(i) \in B_i \quad, \quad 1 \leq i \leq n \quad .
$$
Now form the $0-1$ matrix $P=P_{i,j}$ whose $i$-th row has $1$'s at the columns
belonging to $B_i$, so in this example:
$$
P= \pmatrix{ 
0 & 1 & 1 & 0 & 0 \cr
1 & 1 & 0 & 0 & 1 \cr
0 & 1 & 0 & 0 & 1 \cr
1 & 0 & 0 & 1 & 0 \cr
0 & 1 & 1 & 0 & 0
} 
\quad .
$$

The desired number is simply the {\it permanent} of $P$, $per \,\, P$, in this example
$$
per \,\, P =
per \,\,
 \pmatrix{ 
0 & 1 & 1 & 0 & 0 \cr
1 & 1 & 0 & 0 & 1 \cr
0 & 1 & 0 & 0 & 1 \cr
1 & 0 & 0 & 1 & 0 \cr
0 & 1 & 1 & 0 & 0
} 
=
1 \cdot \,\, per 
 \pmatrix{ 
1 &  0 & 0 & 1 \cr
0 &  0 & 0 & 1 \cr
1 &  0 & 1 & 0 \cr
0 &  1 & 0 & 0
} 
+
1 \cdot \,\, per \,\,
 \pmatrix{ 
1 & 1 &  0 & 1 \cr
0 & 1 &  0 & 1 \cr
1 & 0 & 1 & 0 \cr
0 & 1 & 0 & 0
} = 1+1 =2 \quad .
$$
And indeed, there are two permutations that satisfy the conditions, namely:
$$
\{ 21543 \quad , \quad 31542 \} \quad .
$$

Thanks to Lesley Valiant, we know that for larger $n$, and random (not too small) $B_i$,
this method is {\bf hopeless}.

{\bf A crash course on Rook Polynomials}

Fix, for now, a specific ({\it numeric}), positive integer $n$, and suppose that we are given
any (arbitrary) $n$ subsets $A_i$, $1 \leq i \leq n$, as above, and ask for the number of permutations
of length $n$, $\pi$, such that
$$
\pi(i) \not \in A_i \quad .
$$
So the set of {\it all} guys, AllGuys, is $S_n$, and there are
$$
\sum_{i=1}^{n} |A_i| \quad ,
$$
possible {\it sins} that such a permutation can commit, and we want to enumerate the good guys,
those that have no sins.

There is a natural way of representing permutations in terms of ways of placing $n$ {\it non-attacking}
rooks on an $n \times n$ chessboard. For any permutation $\pi$, let the (sole, of course) rook
on the $i$-th row, be placed at the $\pi(i)$-th column. The condition that $\pi$ is a permutation
guarantees that there is exactly one rook in every column.

Now, given the sets $A_i$ above, cross-out, at each row $i$, the entries of the columns corresponding
to the members of $A_i$, (in other words, you put an $X$ at $P_{ij}$ for $j \in A_i$)
getting a {\it board} , a certain subset of the $n \times n$
chessboard, where it is forbidden to place rooks. For example if, as above
$$
A_1 = \{1,4,5\} \quad , \quad 
A_2 = \{3,4 \} \quad , \quad 
A_3 = \{1,3,4 \} \quad , \quad 
A_4 = \{2,3,5 \} \quad , \quad 
A_5 = \{1,4,5 \} \quad ,
$$
then we have the cross-out board
$$
  \pmatrix{ 
X &  &  & X & X \cr
 &  &  X & X &  \cr
X &  & X & X & \cr
 & X & X &  & X \cr
X & & & X & X
} \quad .
$$

Now it is time to look at the possible `sets of sins', $T$, and the corresponding cardinality of $A_T$.

Let's look at the possible subsets of sins a permutation may be guilty of. 
If the set $T$, of chosen $X$s, has two or more members
on the same row, or two or more members on the same column, then no permutation can be
so sinful, (or else it would not even be a permutation!). On the other hand, if the set, $T$,  of
chosen $X$s has the property that no two of its members are ever on the same row, or on the same column,
and $T$ has  $k$ members, then the cardinality of $A_T$ is $(n-k)!$, since we already know,
{\it for sure}, where $k$ rooks are, and the number of ways of placing the $n-k$ remaining rooks is $(n-k)!$.

This may be viewed as a brand-new combinatorial problem,
that of weight-enumerating the set of placements of {\it any} number of rooks (from $0$ to $n$) on the squares marked by $X$
(the board), in such a way that they are non-attacking,
weighted by $t^k$, where $k$ is the number of rooks participating.
In the above example, we have
$$
R(t)=1+14\,t+63\,{t}^{2}+105\,{t}^{3}+56\,{t}^{4}+6\,{t}^{5} \quad ,
$$
meaning that there is exactly $1$ ways (of course) of placing $0$ rooks, $14$ ways (of course) of placing $1$ rook,
$63$ ways of placing $2$ non-attacking rooks, $\dots$, and $6$ ways of placing $5$ non-attacking rooks.

This is called the {\it rook polynomial} of the board, $R_B(t)$.

We now do the `umbral' substitution $t^i \rightarrow (-1)^i (n-i)!$, or equivalently,
evaluate the integral
$$
\int_{0}^{\infty} e^{-t} t^n R(-\frac{1}{t}) \, dt \quad,
$$
(and get, in this example, with $n=5$), that it equals $2$, confirming the previous
`positive' approach using permanents.

{\bf How to compute Rook Polynomials for a specific (arbitrary) board?}

By recursion of course! It is more convenient to represent a board with a
$0-1$ matrix, where  $1$ denotes $X$ and $0$ denotes a blank square.
So you get the {\it complementary} matrix to the
one above when we did it via permanents. In the above example, 
we have the cross-out board
$$
B=
  \pmatrix{ 
1 & 0 & 0 & 1 & 1 \cr
0 & 0 &  1 & 1 & 0 \cr
1 & 0 & 1 & 1 & 0\cr
0 & 1 & 1 & 0 & 1 \cr
1 & 0 & 0 & 1 & 1
} \quad .
$$

If you input such a $0-1$ matrix, $B$, let's look at the top row. If it is 
{\it not} the all-$0$ row, i.e. there is at least one $1$, let $j$ be the
smallest column with that $1$. Then we have the dilemma `to place or not to place' a rook there.
If we decide {\bf not} to place a rook there, we replace the $1$ by $0$, getting another
$n \times n$ matrix with one fewer $1$, let's call it $B_1$. In this example,
we have the cross-out board
$$
B_1=
  \pmatrix{ 
0 & 0 & 0 & 1 & 1 \cr
0 & 0 &  1 & 1 & 0 \cr
1 & 0 & 1 & 1 & 0\cr
0 & 1 & 1 & 0 & 1 \cr
1 & 0 & 0 & 1 & 1
} \quad .
$$
On the other hand, if we decide to place a rook there, then we must delete the first row, and
the $j$-th column, since no more rooks are allowed there, getting a smaller, $(n-1) \times (n-1)$
matrix $B_2$. In this example, we have
$$
B_2=
 \pmatrix{ 
 0 &  1 & 1 & 0 \cr
 0 & 1 & 1 & 0\cr
 1 & 1 & 0 & 1 \cr
 0 & 0 & 1 & 1
} \quad .
$$
Since $R_B(t)$ is the weight-enumerator, according to the weight $t^{\#Rooks}$, we have the recursion
$$
R_B(t)= R_{B_1}(t) + t R_{B_2(t)} \quad ,
$$
where both $B_1$ and $B_2$ are `smaller' in some sense.

What if the first row only has $0$'s? Then we do the analogous thing to the
first column. What if the first row {\bf and} the first column have all zeros?, then
we delete both the first row and the first column, getting an $(n-1) \times (n-1)$ matrix,
$B'$, and of course $R_B(t)=R_{B'}(t)$.

{\bf But we are mathematicians, NOT accountants, we want GENERAL, Explicit, Symbolic Formulas, Valid for Every $n$}

So we need more structure, and consider `infinite families' of boards.

The grand-daddy of this type of problems is the good-old {\it  Probl\`eme de recontres},
raised way back in 1713 by Pierre R\'emond de Montmort, and solved by Euler,
see the many references in [Sl1].

Here the board is the identity matrix $I_n$, and we have the trivial recurrence
$$
R_n(t)=(t+1)R_{n-1}(t) \quad, \quad R_0(t)=1 \quad ,
$$
leading to the `explicit' expression
$$
R_n(t)=(t+1)^n \quad ,
$$
and hence the number of derangements is
$$
D_n= \int_{0}^{\infty} t^n (1-\frac{1}{t})^n e^{-t} \, = \,
\int_{0}^{\infty} (t-1)^n e^{-t} \quad,
$$
that implies ({\it inter alia} via the not-as-famous-as-it-should-be {\it Almkvist-Zeilberger} 
algorithm[AZ]), the {\it inhomogeneous}  first-order linear recurrence with polynomial coefficients
$$
D_n-nD_{n-1}=(-1)^{n} \quad ,
$$
as well as the {\it homogeneous} linear recurrence with polynomial coefficients
$$
D_n-(n-1) D_{n-1}- (n-1)D_{n-2} =0 \quad .
$$
Both of these recurrences go back, at least, to Euler.

{\bf The Probl\`eme des M\'enages of Edouard Lucas}

A bit more challenging are the boards 
$$
M_1=
\pmatrix{ 
X & X & & & & & & & & \cr
 & X & X & & & & & &   \cr
\dots     \cr
\dots     \cr
 &  &  & & & & & X& X  \cr
 &  &  & & & & & & X} \quad ,
$$
for a {\it straight table}, and 
$$
M_2=
\pmatrix{ 
X & X & & & & & & & & \cr
 & X & X & & & & & & \cr
\dots     \cr
\dots     \cr
 &  &  & & & & & X& X \cr
X &  &  & & & & & & X } \quad ,
$$
for a `round table' (see below).
These come up in the solution of the celebrated {\it  Probl\`eme des M\'enages}
concocted by one of the most {\it amusing} (and deep!) mathematicians of all time,
the inimitable Edouard  Lucas. In his own {\it mots} ([Lu], p. 215, also quoted by Major Percy MacMahon [M], p. 253):

`` {\it Probl\`eme des m\'enages}. -- Des femmes, en nombre $n$, sont range\'es autour d'une table, dans un
ordre d\'etermin\'e ; on demande quel est le nombre des mani\`ers de placer leur maris respectifs, de telle
sorte qu'un homme soit  plac\'e entre deux femmes, sans se trouver
\`a c\^ot\'e de la sienne'' \quad .

In other words, in how many ways can $n$ married couples, sitted around a round table, have
a fun and interesting dinner, where each wife can flirt with both men next to her,
and not sit next to her boring husband, who, being French, probably cheats on her,
and conversely, every husband has opportunities to seduce two other married women.

This problem was independently posed, in a more `serious' context, by Peter Guthrie Tait,
the close collaborator of Lord Kelvin, that somehow came up in his study of knots.

The two sequences that come up (for a straight and round table respectively) are [Sl2] and [Sl3],
see the many references given there. The most satisfactory {\it answer} was given
by {\it monsieur le colonel} Charles Paul Narcisse Moreau, a decorated soldier,
superb amateur mathematician, but not-quite-as-good chess player, and noted
politician and mathematician/engineer Charles Ange Laisant, who derived the
recurrence (for the round table):
$$
A_n = n A_{n-1} + 2 A_{n-2} - (n-4)A_{n-3} - A_{n-4} \quad .
$$

This answer, that only requires $O(N)$ time and memory to compute the
first $N$ terms, is {\bf much better} (at least in the sense of my master, Herb Wilf) than the later `answer' by
Jacques Touchard,
$$
A_n=\sum_{k=0}^{n} (-1)^k \frac{2n}{2n-k} {{2n-k} \choose {k}} (n-k)! \quad ,
$$
that requires $O(N^2)$ time and memory, notwithstanding Gian-Carlo Rota's pronouncement ([Ro], p.340)
that Touchard's `explicit' formula is superior to the recurrences. In fact, using
the {\it Zeilberger algorithm}, one can use
Touchard's `explicit' binomial coefficients sum to give yet-another-proof
of the colonel's and the politician's recurrence.

But why stop here? the great (and fearless!) John Riordan[Ri2], did the next-in-line,
the board
$$
\pmatrix{ 
X & X &X & & & & & & & \cr
 & X & X & X& & & & &   \cr
\dots     \cr
\dots     \cr
 &  &  & & & & X& X& X  \cr
 X &  &  & & & & & X & X \cr
 X& X &  & & & & &  & X
} \quad ,
$$
(and the analogous board for a straight table), that lead to sequences [Sl4] and [Sl5],
but was unable to find  linear recurrences for the sequences themselves, only the generating functions of the
sequences of the rook polynomials, from which one can easily get many terms.

Moving right along, the great enumerator, Earl Glen Whitehead treated[W](by hand!)
the case of {\it four-discordant} permutations, i.e. the board
$$
\pmatrix{ 
X & X &X & X & & & & & & \cr
 & X & X & X& X& & & &   \cr
\dots     \cr
\dots     \cr
 &  &  & & & X & X& X & X  \cr
X &  &  & & &  & X& X & X  \cr
X & X &  & & &  & & X & X  \cr
X& X & X & & & &  &  &  X
} \quad ,
$$
that yielded sequence [Sl6],
once again without finding a recurrence for the sequence itself, but an explicit rational generating function
of the sequence of rook polynomials, and hence as many terms of the enumerating sequence as desired.

[see {\tt http://www.math.rutgers.edu/\~{}zeilberg/tokhniot/oMENAGESegw1f} for a recurrence of order $40$ and
with coefficients that are of degree $2$ in $n$, for the enumerating sequence].

{\bf How did these amazing humans do it?}

They used, either explicitly, or implicitly, the {\it transfer-matrix method},
described in [St], or, other human ingenious, but {\it ad-hoc}, ways, that every time
one had to start from {\it scratch}, and get a {\it Rube Goldberg} monster.

{\bf The General problem}

Since we are not interested in {\it retail} mathematics, only in {\it wholesale}, we want
an {\it algorithm}, {\bf fully implemented}, that can answer it {\it in general}.

{\bf Input}: Any set $S$ of integers (where both negative and positive integers are allowed, as well as $0$, of course).

{\bf Output}: An efficient way to generate the first $1000$ (or whatever) terms of the sequences

(i) $a_S(n)$: the number of $n$-permutations $\pi$ such that $\pi(i) -i \not \in S$

(ii) $b_S(n)$: the number of `circular' $n$-permutations $\pi$ such that $\pi(i) -i \not \in S$,
     where $1$ comes after $n$. For this case, one can assume without loss of generality that
     all the members of $S$ are non-negative, and the smallest member is $0$ (since we are sitting
     around a circular table, you can always move them around).

The best would be a linear recurrence for the enumerating sequence, but failing this (if it is too complicated),
an explicit generating function for the sequence of rook polynomials, from which one can get
as many terms as one desires.

{\bf The Hard Way: Teaching Human Ingenuity to the Computer}

Given such a set $S$ we have to find a {\it scheme} for the rook polynomials of the
$n \times n$ board $B(n)$, implied by $S$. We can call it a `symbolic board'.

Now using the {\it generic} recurrence for specific (numeric) matrices, given above,
one can apply it to the symbolic board $B(n)$, thereby getting other `symbolic boards',
$B'(n)$,  $B''(n)$, that in turn, after you apply the recurrence to them
give rise to {\it yet} more boards $B'''(n)$, $B''''(n)$ etc.

It it is easy to see that if one continues this process, sooner or later, we don't
encounter any new symbolic boards, only previously encountered boards, with $n$ replaced
by $n-1$. This would enable the computer to set-up a system of
linear equations with symbolic coefficients, for the generating functions (w.r.t. $x$) of 
$R_n(t)$,$R'_n(t)$, $R''_n(t)$ etc., let's call them $F(x,t)$, $F'(x,t)$,  $F''(x,t)$, etc.
(see, p. 196 of [Ri1] for a very simple human example, done by hand).

It follows from Cramer's rule that the solutions to this system  of linear equations are rational functions of $(x,t)$.
In particular, our original object of desire, $F(x,t)$ 
$$
 F(x,t) \, = \, \sum_{n=0}^{\infty} R_n(t) x^n  \quad ,
$$
is such a nice  rational function of both $x$ and $t$.
Hence
$$
R_n(t)= \frac{1}{2\pi i} \int_{|x|=1} \frac{F(x,t)}{x^{n+1}} \, dx \quad ,
$$
and hence, the sequence itself: 
$$
A_n=\int_{0}^{\infty} t^n R_n(-\frac{1}{t}) e^{-t} \, dt \quad ,
$$
is given, by the double-integral
$$
A_n=  \frac{1}{2\pi i} \int_{0}^{\infty}  t^n\, e^{-t} \, \int_{|x|=1}  \frac{F(x,-\frac{1}{t})}{x^{n+1}} \, dx \, dt  \quad ,
$$
from which one can (at least in principle, but often also in practice) obtain
a rigorously-derived {\it linear recurrence equation with polynomial coefficients}, using
the multi-Almkvist-Zeilberger algorithm due to Moa Apagodu and Doron Zeilberger[ApZ].

Once a computer has been `taught' how to do `research', i.e. derive the scheme for
the rook polynomials (by introducing dynamically many other auxiliary `boards', but
finitely many of them), and then using the built-in (in Maple) procedure {\tt solve}
that can handle symbolic systems of linear equations, it can do it, at least in principle
for {\it any} finite set $S$ of integers, and thereby reproduce, in a few seconds,
countless hours of human labor (by very smart people,  some of whom were mentioned above),
and then go on to generate new knowledge way beyond the scope of mere humans.

{\bf But}, there is one problem. At this time of writing, one still needs a human to do the
teaching! In other words, design an algorithm, and implement it, that
does the above for any set $S$. And, don't expect too much!
As the sets $S$ get bigger, even computers would refuse to do it!
Computers are only a few orders of magnitude better than humans.
Programming the computer to do such research requires at least
as much ingenuity, (and hard work!) as doing special cases by hand. 
But, if you are like me and are lazy, and don't feel like wasting weeks writing
a long  computer program, you can use the {\it Zeilberger Gordian knot}.
Replace {\it ingenuity} by {\it meta-ingenuity}.

{\bf The `Gordian Knot' Way: Long Live the Empirical (yet rigorous) approach!}

Since we know, {\it a priori}, by general `hand-waving' (that can be easily made rigorous), that
there {\bf exists} a scheme, or in the language of Richard Stanley [St], there {\it exists}
a transfer matrix, we know {\it a priori}, that the sequence of rook polynomials for $S$,
$R_n(t)$, satisfies {\it some} linear recurrence equation with {\it constant} coefficients,
i.e. belongs to the $C$-finite ansatz [Z][KP].
This means that the (ordinary) generating function 
with respect to $x$ is a {\it rational function} of $x$ and $t$. 
Hence we can use the {\it straightforward} specific procedure, mentioned above, for finding $R_B(t)$
for {\it any} board, to crank out the first few terms of the sequence
of rook polynomials for the $n \times n$ boards implied by $S$,
and then {\it guess} the recurrence, that we {\bf know} exists!
Then the computer can either find the implied rational generating function, and
use the multi-Almkvist-Zeilberger algorithm, but a much better way, once the
linear recurrence equation for the sequence of rook polynomials is known, is to crank out many more terms, and
then use them (by the umbral operator $t^i \rightarrow (-1)^i(n-i)!$) to
crank out sufficiently many terms of the sequence $a_S(n)$  (or $b_S(n)$), and
{\it guess} a linear recurrence equation with {\it polynomial} coefficients
(i.e. use the {\it holonomic ansatz}), that once again can be justified fully rigorously (if desired),
since we know that it {\bf exists}, and one can easily find {\it a priori} upper bounds for  the order and degree.

{\bf The Maple package MENAGES}

All this (and much more!) is implemented in the Maple package {\tt MENAGES} available directly from

{\tt http://www.math.rutgers.edu/\~{}zeilberg/tokhniot/MENAGES} \quad ,

or by clicking on the indicated link in the front of the present article

{\tt http://www.math.rutgers.edu/\~{}zeilberg/mamarim/mamarimhtml/menages.html } \quad ,

that contains lots of input and output files. To access the on-line help just
type {\tt ezra(); }. I will only mention a few of the available procedures.

$\bullet$ {\tt RP(M,t);}: inputs  a square $0-1$ matrix, M, and a variable $t$, and outputs the
rook polynomial of $M$, where the $1$'s denote the $X$'s and the $0$'s the other squares of the
$n \times n$ board.

$\bullet$ {\tt RookPrec(S,t);}: inputs a set of integers $S$ and a variable $t$ and outputs the linear recurrence
equation with constant coefficients satisfied by the rook polynomials of the board implied by the
problem of counting $n$-permutations $\pi$ such that $\pi(i) - i \not \in S$.

For example,  `{\tt RookPrec($\{0,1\}$ , t);}' outputs  $[[1+t,1+3\,t+{t}^{2}],[1+2\,t,-{t}^{2}]]$, which is
the package's way of telling you that the sequence of Rook polynomials, $R_n(t)$ for an $n \times n$
board, satisfies the recurrence
$$
R_n(t)=(1+2t)R_{n-1}(t) -t^2 R_{n-2}(t) \quad , \quad
$$
subject to the initial conditions
$$
R_1(t)=1+t \quad , \quad R_2(t)= 1+ 3t+t^2 \quad .
$$

$\bullet$ {\tt RookPcrec(S,t);}: is the analogous procedure for {\it circular} permutations (diners around a round table)

$\bullet$ {\tt Seq(S,N);}: inputs a set of integers $S$ and outputs the
the first $N$ terms of the sequence enumerating $n$-permutations such that
$\pi(i)-i \not \in S$, starting at $n=1$. For example

{\tt Seq($\{0,1\}$,21);}

yields  [Sl2] (i.e. {\tt http://oeis.org/A000271}).

$\bullet$ {\tt SeqC(S,N);}: inputs a set of integers $S$ and outputs the
the first $N$ terms of the sequence enumerating $n$-circular permutations such that
$\pi(i)-i \not \in S$ (where $1$ follows $n$) starting at $n=1$. For example

{\tt SeqC($\{0,1\}$,21);}

yields  [Sl3] (i.e. {\tt http://oeis.org/A000179}), and all these terms, except the last one are
already in Lucas' classic[Lu](p. 495).

$\bullet$ {\tt InfoE(S,n,N,MaxC,K,L1,L2);}: inputs a set of  integers {\tt S}, 
symbols $n$ and $N$ (where $N$ is the shift operator in $n$)
and outputs a list consisting of

(i) A list of  {\tt L1} terms whose $n$-th term is the number of permutations $\pi$ of $\{1, ...,n\}$ such that
$\pi(i)-i \not \in S$, let's call it $a(n)$

(ii) a linear recurrence operator of complexity {\tt ORDER+DEGREE} $\leq$ {\tt MaxC}, 
of minimal order, $ORDER$, annihilating the sequence, or FAIL, if none exists.

(iii) {\tt a(L2)} (if the second component is not FAIL), otherwise FAIL.

(iv) The asymptotic expression for $a(n)/n!$ (or FAIL), to order {\tt K}.

$\bullet$ {\tt InfoEC(S,n,N,MaxC,K,L1,L2);}: ditto for {\bf circular} permutations, except that w.l.o.g
the set $S$ can consist of non-negative integers and  its smallest element is $0$.

{\tt InfoEV},  {\tt InfoECV} are verbose versions, and procedures {\tt SeferE} and {\tt SeferEC} output
webbooks. See the on-line help.

Other procedures that output articles are   {\tt Mamar}, {\tt MamarC}, and {\tt MamarDiscordant}. Try them out!

{\bf The Maple package BALTIC}

The complementary, easier (at least conceptually), problem of finding generating functions for
the enumeration of permutations $\pi$, where $\pi(i)-i \in S$ for a prescribed set of integers $S$, is
treated in the small Maple package {\tt BALTIC}, reproducing empirically-yet-rigorously
results of Lehmer[Le], Baltic[B], and Stanley [St] (pp. 252-253). In particular 

{\tt GFbaltic($\{$-2,-1,1,2 $\}$,t,20); }

finds in  $0.024$ seconds, that the generating function for the enumerating sequence for permutations
$\pi$ with $\pi(i) -i =$ $\pm 1$ or $\pm 2$ is
$$
{\frac {1-t}{1-t-{t}^{2}-{t}^{3}-{t}^{4}+{t}^{5}}} \quad
$$
(given at the bottom of page 252 of [St]), that took Richard Stanley, using  admittedly elegant
human reasoning, quite a bit longer, and

{\tt GFbaltic($\{$-2,-1,0, 1,2 $\}$,t,20); } yields
$$
\frac{1-t}{1-2\,t-2\,{t}^{3}+{t}^{5}} \quad ,
$$
given at the top of p. 253 of [St]. For far deeper generating functions see:

{\tt http://www.math.rutgers.edu/\~{}zeilberg/tokhniot/oBALTIC2b} \quad .

{\bf Conclusion: It is Time to Stop to Think and  Start to Meta-Think}

Three of the greatest combinatorics gurus of the 20th century , John Riordan, Gian-Carlo Rota, and Richard Stanley 
([Ri1][Ro][St]) taught us how to apply  the {\it transfer-matrix method}, and the {\it principle of Inclusion-Exclusion},
by using {\it insight} and {\it ingenuity} ,
to solve combinatorial problems. Alas, for each separate problem, one
needs {\it both} {\bf inspiration} and {\bf perspiration}, wasting
lots of very-smart-people's precious time.  In this paper we did a bit of
{\bf meta-thinking}, and realized, by `hand-waving' (but fully rigorously!),
that {\bf because} each of these problems may be used by
combining the principle of inclusion and exclusion (via rook polynomials) and the transfer-matrix method,
it follows {\bf a priori}, that the sequence of rook polynomials,
for any `forbidden' set $S$, belongs to
the $C$-finite ansatz, i.e., satisfies {\it some} linear recurrence equation
with {\it constant} coefficients. Having realized that, we can ask our
beloved silicon collaborators to crank out the first few terms,
and using straightforward guessing, {\it guess} the recurrence, that
immediately becomes {\it proved} because of general principles. Then it follows
from `general holonomic nonsense' that the enumerating sequences themselves are
{\it holonomic} (aka $P$-recursive), and once again, a recurrence can be guessed,
that {\it a posteriori} is fully rigorous (not that I care!).

The moral is: an ounce of meta-thinking (plus help from our computers) is worth a pound of (human) thinking!

{\bf Appendix: Some Recommended Webbooks by Shalosh B. Ekhad generated by MENAGES}

$\bullet$ To see the generating functions for rook polynomials for $r$-discordant 
permutations for $r$ from $1$ to $4$, 
reproducing in $12$ seconds, the labor of Euler ($r=1$), Lucas (and Laisant, Moreau, Touchard, Kaplasnky and many other smart people) 
($r=2$), Riordan[Ri2] ($r=3$ plus we got a brand-new recurrence for the enumerating sequence itself!), and Whitehead[W]($r=4$),
look at:

{\tt http://www.math.rutgers.edu/\~{}zeilberg/tokhniot/oMENAGES7} \quad .

$\bullet$ To see the generating functions for rook polynomials for $r$-discordant 
permutations for $r=5,6$ (in addition to the above $1 \leq r \leq 4$), look at:

{\tt http://www.math.rutgers.edu/\~{}zeilberg/tokhniot/oMENAGES7a} \quad .

$\bullet$ To see a webbook that gives you the generating functions for rook polynomials for enumerating permutations 
$\pi$ of $\{0,1,..,n-1\}$ such that $\pi(i)-i \bmod n$ is never in the set $S$, 
for all subsets $S$ of $\{0,1,2,3,4\}$ that include $0$, and in many cases, nice recurrences for the enumerating sequences themselves,
look at:

{\tt http://www.math.rutgers.edu/\~{}zeilberg/tokhniot/oMENAGES8} \quad .

$\bullet$ To see a webbook that gives you the generating functions for rook polynomials for enumerating permutations 
$\pi$ of $\{1,..,n\}$ such that $\pi(i)-i$ is never in the set $S$, 
for all subsets $S$ of $\{-2,-1,0,1,2\}$ of cardinality at least $2$ (some of the theorems are
trivially equivalent, by symmetry, but who cares?), look at:

{\tt http://www.math.rutgers.edu/\~{}zeilberg/tokhniot/oMENAGES9}  \quad .

{\bf References}

[AlZ] G. Almkvist and D. Zeilberger,
{\it The method of differentiating under the integral sign}, 
J. Symbolic Computation {\bf 10} (1990), 571-591.
\hfill\break
{\tt http://www.math.rutgers.edu/\~{}zeilberg/mamarim/mamarimPDF/duis.pdf} \quad .

[ApZ] Moa Apagodu and D. Zeilberger, {\it Multi-Variable Zeilberger and Almkvist-Zeilberger Algorithms and the Sharpening of Wilf-Zeilberger Theory},
Adv. Appl. Math. {\bf 37}(2006), (Special issue in honor of Amitai Regev), 139-152.
\hfill\break
{\tt http://www.math.rutgers.edu/\~{}zeilberg/mamarim/mamarimhtml/multiZ.html} \quad .

[B] Vladimir Baltic, 
{\it On the number of certain types of strongly restricted permutations}, 
Applicable Analysis and Discrete Mathematics {\bf 4} (2010), 119-135.
\hfill\break
{\tt http://pefmath.etf.rs/vol4num1/AADM-Vol4-No1-119-135.pdf } \quad .

[KP] Manuel Kauers and Peter Paule, ``{\it The Concrete Tetrahedron}'', Springer, 2011 \quad .

[Le] D. H. Lehmer, {\it Permutations with strongly restricted displacements},
                   in ``Combinatorial theory and its applications, II (Proc. Colloq., Balatonfured, 1969)'', 755-770. North-Holland, Amsterdam, 1970.

[Lu] Edouard Lucas, {\it ``Th\'eorie Des Nombres''},
Gauthier-Villars et Fils, Paris, 1891.

[M] Major Percy A. MacMahon, F.R.S., D.Sc., LL.D., {\it `` Combinatory Analysis''}, volume I,
Cambridge University Press, 1917. Reprinted (together with volume II), by Chelsea, 1984.

[Ri1] John Riordan, {\it ``Introduction to Combinatorial Analysis''}, John Wiley, 1958.
Reprinted by Dover, 2002.

[Ri2] John Riordan, {\it The enumeration of permutations with three-ply staircase restrictions},
unpublished memorandum, Bell Telephone Laboratories, Murray Hill, NJ, Oct 1963.
\hfill\break
Available (thanks to Neil Sloane!) from {\tt http://oeis.org/A001883/a001883\_21.pdf} \quad .

[Ro] Gian-Carlo Rota, {\it On the Foundation of Combinatorial Theory I. Theory of  M\"obius Functions},
Zeitschrift  Wahrscheinlichkeitstheorie {\bf 2}(1964), 340-368, reprinted in
``Classic Papers in Combinatorics'', (Ira Gessel and Gian-Carlo Rota, editors), Birkhauser, 1987.

[Sl1] Neil Sloane, ``{\it The On-Line Encyclopedia of Integer Sequences (OEIS)}'',\hfill\break
Sequence {\tt http://oeis.org/A000166} \quad .

[Sl2] Neil Sloane, ``{\it The On-Line Encyclopedia of Integer Sequences (OEIS)}'',\hfill\break
Sequence {\tt http://oeis.org/A000271 } \quad .

[Sl3] Neil Sloane, ``{\it The On-Line Encyclopedia of Integer Sequences (OEIS)}'',\hfill\break
Sequence {\tt http://oeis.org/A000179 } \quad .

[Sl4] Neil Sloane, ``{\it The On-Line Encyclopedia of Integer Sequences (OEIS)}'',\hfill\break
Sequence {\tt http://oeis.org/A000183 } \quad .

[Sl5] Neil Sloane, ``{\it The On-Line Encyclopedia of Integer Sequences (OEIS)}'',\hfill\break
Sequence {\tt http://oeis.org/A001887 } \quad .

[Sl6] Neil Sloane, ``{\it The On-Line Encyclopedia of Integer Sequences (OEIS)}'',\hfill\break
Sequence {\tt http://oeis.org/A004307 } \quad .

[St] Richard Stanely, ``{\it Enumerative Combinatorics, volume I}'' (first edition), Wadsworth \& Brooks/Cole, 1986.

[W] Earl Glen Whitehead, Jr., {\it Four-discordant permutations}, J. Austral. Math. Soc. Ser. A {\bf 28}(1979), 369-377.

[Z] Doron Zeilberger, {\it The C-finite ansatz}, Ramanujan Journal {\bf 31} (2013), 23-32.
\hfill\break
{\tt http://www.math.rutgers.edu/\~{}zeilberg/mamarim/mamarimhtml/cfinite.html}  \quad .
\end